\documentclass[11pt, a4paper]{amsart}

\usepackage[T1]{fontenc}
\usepackage{ae}

\newcommand\OO{\mathcal{O}}
\newcommand\QQ{\mathbb{Q}}
\newcommand\ZZ{\mathbb{Z}}
\newcommand\PP{\mathbb{P}}
\newcommand\RR{\mathbb{R}}
\newcommand\Pone{{\PP^1}}
\newcommand\Ptwo{{\PP^2}}
\newcommand\xx{\mathbf{x}}
\newcommand\yy{\mathbf{y}}
\newcommand\dd{\mathrm{d}}

\newcommand{\tS}{{\widetilde S}}
\newcommand{\Aone}{{\mathbf A}_1}
\newcommand{\Atwo}{{\mathbf A}_2}
\newcommand{\Athree}{{\mathbf A}_3}
\newcommand{\Afour}{{\mathbf A}_4}
\newcommand{\Afive}{{\mathbf A}_5}
\newcommand{\Dfour}{{\mathbf D}_4}
\newcommand{\Dfive}{{\mathbf D}_5}
\newcommand{\Esix}{{\mathbf E}_6}
\newcommand{\Eseven}{{\mathbf E}_7}
\newcommand{\Eeight}{{\mathbf E}_8}
\DeclareMathOperator{\Pic}{Pic}
\newcommand\Ceff{{\Lambda_\mathrm{eff}}}
\newcommand\Cdual{{\Lambda_\mathrm{eff}^\vee}}
\DeclareMathOperator{\vol}{Vol}
\newcommand\polymake{\texttt{Polymake}}

\newtheorem{theorem}{Theorem}
\newtheorem{lemma}[theorem]{Lemma}
\newtheorem{question}[theorem]{Question}

\theoremstyle{definition}
\newtheorem{definition}[theorem]{Definition} 

\theoremstyle{remark}
\newtheorem{remark}[theorem]{Remark}

\begin{document}

\title[A constant arising in Manin's conjecture for Del Pezzo surfaces]
{On a constant arising in Manin's conjecture for Del Pezzo surfaces}

\author{Ulrich Derenthal} 

\address{Institut f\"ur Mathematik, Universit\"at
  Z\"urich, Winterthurerstrasse 190, 8057 Z\"urich, Switzerland}

\email{ulrich.derenthal@math.unizh.ch}


%14J26: Rational and ruled surfaces
%14G05: Rational points
%52B05: Polytopes and polyhedra - Combinatorial properties 
%(number of faces, shortest paths, etc.)

\begin{abstract}
  For split smooth Del Pezzo surfaces, we analyse the structure of the
  effective cone and prove a recursive formula for the value of
  $\alpha$, appearing in the leading constant as predicted by Peyre of
  Manin's conjecture on the number of rational points of bounded
  height on the surface.  Furthermore, we calculate $\alpha$ for all
  singular Del Pezzo surfaces of degree $\ge 3$.
\end{abstract}

\maketitle

\tableofcontents

\section{Introduction}

Over the field $\QQ$ of rational numbers, a split smooth Del Pezzo
surface $S$ is $\Ptwo$, $\Pone \times \Pone$, or, for $1 \le r \le 8$,
the blow-up $S_r$ of $\Ptwo$ in $r$ points which are defined over
$\QQ$ and are in \emph{general}\footnote{I.e., no three points on a
  line, no six on a conic, no eight on a cubic curve with one of them
  a singularity of that curve.} position. For the degree of
$S$, we have $\deg(\Ptwo) = 9$, $\deg(\Pone \times \Pone) = 8$, and
$\deg(S_r) = 9-r$.

An important object associated to $S$ is the
\emph{effective cone} $\Ceff(S)$, i.e., the convex cone in
\[\Pic(S)_\RR := \Pic(S) \otimes_\ZZ \RR\] which is generated by the
classes of effective divisors. Note that we identify divisors and their
classes in $\Pic(S)$ when this cannot lead to confusion.

On $\Pic(S) \cong \ZZ^{10-\deg(S)}$, we have the \emph{intersection
  form} $(\cdot, \cdot)$, which is a non-degenerate symmetric bilinear
form. The \emph{anticanonical map} $\phi: S \to \PP^{\deg(S)}$ is
given by the \emph{anticanonical class} $-K_S =
\phi^*(\OO_{\PP^{\deg(S)}}(1))$. It is an embedding for $\deg(S) \ge
3$.

For $S = S_r$, the anticanonical class is \[-K_r := -K_{S_r} =
3H-(E_1+\dots+E_r),\] where $H$ is the transform of a general line in
$\Ptwo$, and $E_1, \dots, E_r$ are the exceptional divisors obtained
by blowing up the $r$ points in $\Ptwo$. By
\cite[Corollary~3.3]{MR2029863}, for $r \ge 2$, $\Ceff(S_r)$ is
generated by the \emph{$(-1)$-curves}, i.e., prime divisors $D$ whose
self-intersection number $(D,D)$ is $-1$.  For $r \le 6$, the
anticanonical embedding $\phi_r = \phi$ maps the $(-1)$-curves on
$S_r$ exactly to the lines on $\phi_r(S_r) \subset \PP^{9-r}$. Note
that $H, E_1, \dots, E_r$ give a basis of $\Pic(S_r) \cong \ZZ^{r+1}$.

Starting in 1989, Manin initiated a program to study the number of
rational points on certain varieties which can be stated in case of a
split smooth Del Pezzo surface $S$ of degree $d$ as follows:

\begin{question}
  Let $U$ be the complement of the $(-1)$-curves in $S$, and let
  \[H: S(\QQ) \to \ZZ_{>0}\] be the anticanonical height, i.e., $H(\xx) =
  \max\{|x_0|, \dots, |x_d|\}$ where the image in $\PP^d(\QQ)$ of $\xx \in
S(\QQ)$ under the anticanonical map $\phi: S \to \PP^d$ is represented by
integral and coprime coordinates $(x_0, \dots, x_d)$.

  What is the asymptotic behavior of the \emph{number of rational
    points of bounded height} \[N_{U,H}(B) := \#\{\xx \in U(\QQ) \mid
  H(\xx) \le B\}\] as $B \to \infty$?
\end{question}

By Manin's conjecture \cite{MR89m:11060}, which was formulated for the
much larger class of Fano varieties, the following answer is expected:
\[N_{U,H}(B) \sim c_{S,H} \cdot B \cdot (\log B)^{9-d}.\] The leading
constant $c_{S,H}$ has received a conjectural interpretation by
Peyre \cite{MR1340296}: it is expected to be the product \[c_{S,H} =
\alpha(S)\cdot \beta(S) \cdot \omega_H(S),\] where $\alpha(S)$
is a constant related to the geometry of $S$, while $\beta(S)$ is
a cohomological constant which is always $1$ for split Del Pezzo
surfaces, and $\omega_H(S)$ is related to the densities of rational
points on $S$ over $\RR$ and modulo $p$ for all primes $p$.

So far, Manin's conjecture has been proved for split smooth Del Pezzo
surfaces in the cases $d \ge 6$ in the context of the more general
proof for toric varieties \cite{MR1620682}, and for a specific surface
in the case $d=5$ \cite{MR1909606}.

\medskip

The purpose of this note is to look more closely at the constant
$\alpha(S)$. Its definition is due to Peyre
(\cite[Definition~2.4]{MR1340296}; see \cite[Section~6]{MR1875177} for
more details):

\begin{definition}\label{def:alpha}
  Let $\Ceff(S)$ be the effective cone, $\Cdual(S)$ its dual cone
  (with respect to the intersection form) of \emph{nef} divisor
  classes, and $-K_S$ the anticanonical class on $S$. Then we define
  \[\alpha(S) := \vol(P(S)),\]
  where \[P(S):=\{x \in \Cdual(S) \mid (-K_S, x) = 1\}\] is a polytope whose
  volume is calculated using the Lebesgue measure on the hyperplane $\{x \in
  \Pic(S)_\RR^\vee \mid (-K_S, x) = 1\}$ which is defined by the $(9-d)$-form
  $\dd \xx$ such that $\dd \xx \wedge \dd \omega = \dd \yy$, where $\dd \yy$
  is the form corresponding to the natural Lebesgue measure on
  $\Pic(S)^\vee_\RR$ and $\dd\omega$ is the linear form defined by $-K_S$ on
  $\Pic(S)^\vee_\RR$.
\end{definition}

For large $d$, the calculation of $\alpha(S)$ can be carried out
directly by hand (see \cite[Section~1.3]{MR1909606} for the case
$d=5$). For small $d$, especially for $S_8$ of degree $d=1$, a direct
calculation seems to be currently impossible even with the help of
software like $\polymake$ \cite{MR1785292}. In this case, the cone
$\Ceff(S_8)$ has 240 generators, while $\Cdual(S_8)$ has 19440
generators. A direct calculation of $\alpha(S_8)$ would require a
triangulation of $\Cdual(S_8)$, which seems to be out of reach for
today's software and hardware.

Therefore, we need a more detailed knowledge of $\Cdual(S)$. For $S =
S_r$ and $r \ge 3$, we have an action of a Weyl group $W_r$ on
$\Pic(S_r)$; see Table~\ref{tab:smooth} for the type of $W_r$ and
\cite[Section~2]{MR2029863} for details. Our main result which will
allow us to compute $\alpha(S_r)$ recursively is:

\begin{theorem}\label{thm:cone}
  Let $r \ge 3$.  The nef cone $\Cdual(S_r)$ has $N_r$ faces, where
  $N_r$ is the number of $(-1)$-curves on $S_r$. Each face is
  isomorphic to $\Cdual(S_{r-1})$. The Weyl group $W_r$ acts
  transitively on the faces and leaves $-K_r$ in the interior of
  $\Cdual(S_r)$ fixed.
\end{theorem}

This observation is a crucial step in the proof of the following
recursive formula for $\alpha(S_r)$; see Table~\ref{tab:smooth} for
the values of $\alpha(S_r)$ and $N_r$:

\begin{theorem}\label{thm:main}
  Let $S_r$ be a split smooth Del Pezzo surface of degree $9-r$ which
  is the blow-up of $\Ptwo$ in $r$ points in general position. Let
  $N_r$ be the number of $(-1)$-curves on $S_r$. We have $\alpha(S_2) =
  1/24$ and\[\alpha(S_r) = \frac{N_r\cdot \alpha(S_{r-1})}{r \cdot
    (9-r)}\] for $3 \le r \le 8$.  Furthermore, $\alpha(S_1) = 1/6$,
  $\alpha(\Pone \times \Pone) = 1/4$, and $\alpha(\Ptwo) = 1$.
\end{theorem}

\begin{table}[ht]
  \centering
  \[\begin{array}{|c||c|c|c|c|c|c|c|}
    \hline
    r & 2 & 3 & 4 & 5 & 6 & 7 & 8\\
    \hline\hline
    \text{type of $W_r$} & & \Atwo \times \Aone & \Afour & \Dfive & \Esix & \Eseven & \Eeight\\ 
    \hline
    N_r & 3 & 6 & 10 & 16 & 27 & 56 & 240\\
    \hline
    \alpha(S_r) & 1/24 & 1/72 & 1/144 & 1/180 & 1/120 & 1/30 & 1\\
    \hline
  \end{array}\]
  \caption{Smooth Del Pezzo surfaces}
  \label{tab:smooth}
\end{table}

Next, we consider split singular Del Pezzo surfaces $S$ whose
singularities are rational double points. Besides the case where $S$
is the Hirzebruch surface $F_2$ of degree $8$, their \emph{minimal
  desingularizations} $\tS$ are obtained as follows: we perform a
series of $r \le 8$ blow-ups of $\Ptwo$ resulting in $\tS = \tS_r$,
where in at least one step, we blow up a point on a $(-1)$-curve
(resulting in $(-2)$-curves, i.e., prime divisors with self intersection
number $-2$), where
the only restriction for the choice of the blown-up point is that we
never blow up a point on a $(-2)$-curve (therefore, no prime divisors with
self intersection number smaller than $-2$ can occur). Contracting the
$(-2)$-curves on $\tS_r$ results in $S=S_r$ of degree $9-r$.

With some minor modifications in its formulation (mostly replacing $S$
by $\tS$ where appropriate), Manin's conjecture is expected to hold
for singular Del Pezzo surfaces as well; see \cite{math.NT/0511041}
and the introduction of \cite{math.NT/0604193} for details and an
overview of the current progress.  The definition of $\alpha(S)$ is
also almost the same: we must consider $-K_\tS$, $\Ceff(\tS)$, and
$\Cdual(\tS)$ in $\Pic(\tS)$ of rank $10-\deg(S)$.  Note that
$\Ceff(\tS)$ is generated by the \emph{negative curves} (i.e., the
$(-1)$- and $(-2)$-curves) in the singular case if $\deg(S) \le 7$.
The value of $\alpha(S)$ depends not only on the
degree of $S$, but also on the type of singularities on $S$.

In Section~\ref{sec:singular}, we list $\alpha(S)$ for each type of
singular Del Pezzo surface of degree at least 3; see
Tables~\ref{tab:degree_7}, \ref{tab:degree_6}, \ref{tab:degree_5},
\ref{tab:degree_4}, and \ref{tab:degree_3}. For most types, the
calculation was performed with the help of the data given in
\cite{math.AG/0604194} and \cite[Chapter~6]{thesis}.

For some examples of the calculation of $\alpha(S_r)$ for non-split
Del Pezzo surfaces, see \cite{MR1340296}, \cite[Section~6]{MR1875177},
\cite{MR2099200}, and \cite{math.NT/0502510}.

\medskip

\noindent{\textbf{Acknowledgments.}}
I am grateful to T.\ D.\ Browning for several useful comments.

\section{Smooth Del Pezzo surfaces}\label{sec:smooth}

Let $S_r$ be the blow-up of $\Ptwo$ in $r$ points in general position.

\begin{lemma}\label{lem:effective_cone}
  Let $2 \le r \le 8$. The effective cone $\Ceff(S_r)$ is generated
  over $\RR$ by the $(-1)$-curves, whose number is $N_r$ as listed in
  Table~\ref{tab:smooth}.
\end{lemma}

\begin{proof}
  See \cite[Corollary~3.3]{MR2029863}. Their number $N_r$ can be found
  in \cite[Theorem~2.1]{MR2029863}. (Note that, for $r=8$, the
  semigroup of classes of effective divisors is generated by the
  $(-1)$-curves together with $-K_8$.)
\end{proof}

\begin{lemma}\label{lem:negative_curves}
  Let $E$ be a $(-1)$-curve on $S_r$ for $r \ge 3$. If $D \in \Pic(S_r)$
  fulfills $(D,E) = 0$ and $(D,E') \ge 0$ for all $(-1)$-curves $E'$
  such that $(E,E') = 0$, then $D$ is nef.
\end{lemma}

\begin{proof}
  As the $(-1)$-curves generate the effective cone, we must show that
  $(D,E') \ge 0$ also holds for all $(-1)$-curves $E'$, regardless of
  the value of $(E,E')$.
  
  If $(E,E') < 0$, then $E'=E$, and $(D,E) = 0$. If $(E,E') = 0$, then
  $(D,E') \ge 0$ by assumption. 
  
  We proceed by induction on $n = (E,E')$. If $n = 1$, then $E+E'$ is a ruling
  as in \cite[Definition~4.6]{MR2029863}. The case $n = 2$ occurs for $r \in
  \{7,8\}$, and $n = 3$ is possible only for $r = 8$; furthermore, $n \ge 4$
  is impossible. The divisor $E+E'$ can be written in at least two ways as the
  sum of two negative curves (see \cite[Section~4]{MR2029863} for rulings, and
  \cite[Sections~3.4,~3.5]{thesis} for $n \in \{2,3\}$), say $E+E' =
  E_1+E_2$, where $E \notin \{E_1, E_2\}$. Then
  \[(E,E_1)+(E,E_2) = (E,E')+(E,E) = n-1,\] where $(E,E_1)$ and
  $(E,E_2)$ are both non-negative. Therefore, the induction hypothesis
  holds for $E_1, E_2$, and \[(D,E') = (D, E+E') = (D, E_1+E_2) =
  (D,E_1) + (D,E_2) \ge 0\] completes the induction.
\end{proof}

For any $D \in \Pic(S_r)$ and $c \in \RR$, we define
\[D^{=c} := \{D' \in \Pic(S_r)_\RR \mid (D, D') = c\},\]
and similarly,
\[D^{\ge c} := \{D' \in \Pic(S_r)_\RR \mid (D, D') \ge c\}.\]

\begin{proof}[Proof of Theorem \ref{thm:cone}]
  By definition, $\Cdual(S_r)$ is the intersection of the half spaces $E^{\ge
    0}$ for all generators $E$ of $\Ceff(S_r)$, which are exactly the
  $(-1)$-curves by Lemma~\ref{lem:effective_cone}. By
  \cite[Lemma~5.3]{MR1941576}, $W_r$ acts transitively on the $(-1)$-curves.
  This symmetry implies that each $(-1)$-curve $E$ defines a proper face
  $F_E:=\Ceff(S_r) \cap E^{=0}$, and that $W_r$ acts transitively on the set
  of faces $\{F_E \mid \text{$E$ is a $(-1)$-curve}\}$.

Consider $S_r$ as the blow-up of $S_{r-1}$ in one point, resulting in
the exceptional divisor $E_r$. Then \[\Pic(S_r) = \Pic(S_{r-1}) \oplus
\ZZ \cdot E_r\] is an orthogonal sum.

We claim that $F_{E_r} = \Ceff(S_{r-1})$, where we regard
$\Ceff(S_{r-1}) \subset \Pic(S_{r-1})$ as embedded into $\Pic(S_r)$.

Indeed, if $D \in \Ceff(S_{r-1})$, then $(D,E_r) = 0$, and $(D, E) \ge
0$ for all $(-1)$-curves $E$ of $S_{r-1}$, which are exactly the
$(-1)$-curves of $S_r$ with $(E,E_r) = 0$. By
Lemma~\ref{lem:negative_curves}, we have $(D,E) \ge 0$ for all
$(-1)$-curves of $S_r$.

On the other hand, if $D \in E_r^{=0}$, then $D \in \Pic(S_{r-1})$.
If $D \in \Ceff(S_r)$, then $(D, E) \ge 0$ for all $(-1)$-curves of
$S_r$, which includes the $(-1)$-curves of $S_{r-1}$, proving the
other direction.

The root system corresponding to $W_r$ is \[R_r = \{D \in \Pic(S_r)
\mid (D,D) = -2, (D, -K_r) = 0\}.\] Since $W_r$ is generated by the
reflections $E \mapsto E+(D,E) \cdot D$ corresponding to the roots $D
\in R_r$, the anticanonical class $-K_r$ is fixed under $W_r$. This
completes the proof of Theorem~\ref{thm:cone}.
\end{proof}  

\begin{proof}[Proof of Theorem \ref{thm:main}]
The polytope $P_r:=P(S_r)$ whose volume is $\alpha(S_r)$ (see
Definition~\ref{def:alpha}) is the intersection of the $N_r$
half-spaces $E^{\ge 0}$ (where $E$ runs through the $(-1)$-curves of
$S_r$) in the $r$-dimensional space $-K_r^{=1}$.
  
Note that $(-K_r,-K_r) = 9-r$. Therefore, $Q:=\frac 1{9-r}\cdot (-K_r) \in
-K_r^{=1}$, and since $(-K_r, E) = 1$ for any $(-1)$-curve $E$, the
point $Q$ is in the interior of $P_r$.
  
Consider the convex hull $P_E$ of $Q$ and the face $P_r \cap E^{=0}$
of $P_r$ corresponding to $E$. Then $P_r$ is the union of the $P_E$
for all $(-1)$-curves $E$, and since their intersections are
lower-dimensional, \[\vol(P_r) = \sum_E \vol(P_E).\]
  
As the intersection form and $-K_r$ are invariant under the Weyl group
$W_r$, it acts on $-K_r^{=1}$ and therefore on $P_r$. As in
Theorem~\ref{thm:cone}, it permutes the faces of $P_r$ transitively.
As $Q$ is fixed under $W_r$ and the volume is invariant under $W_r$,
we have $\vol(P_r) = N_r \cdot \vol(P_E)$ for any $(-1)$-curve $E$.
  
As in the proof of Theorem~\ref{thm:cone}, we consider $S_r$ as the
blow-up of $S_{r-1}$ in one point, resulting in the exceptional
divisor $E_r$, with the orthogonal sum \[\Pic(S_r) = \Pic(S_{r-1})
\oplus \ZZ \cdot E_r.\]
  
We claim that $P_r \cap E_r^{=0} = P_{r-1}$. In view of
Theorem~\ref{thm:cone}, it remains to prove that $(D,-K_{r-1}) = 1$ is
equivalent to $(D,-K_r) = 1$ on $E_r^{=0}$. This follows directly from
$-K_r = -K_{r-1} - E_r$.
  
Therefore, $P_{E_r}$ is a cone over the $(r-1)$-dimensional polytope
$P_{r-1}$ in the $r$-dimensional space $-K_r^{=1}$. A cone of height
$1$ over $P_{r-1}$ has volume $\vol(P_{r-1}) / r$. As $E_r$ is
orthonormal to $\Pic(S_{r-1})$, and $(-K_r, E_r) = 1$, the distance of
$Q$ to $P_{r-1}$ is $1/(9-r)$. Therefore,
\[\vol(P_{E_r}) = \frac{\vol(P_{r-1})}{r\cdot(9-r)}.\]
  
Together with $\alpha(S_{r-1}) = \vol(P_{r-1})$ and $\vol(P_r) = N_r
\cdot \vol(P_{E_r})$, this completes the proof of the recursive
formula.
  
For $r=2$, we have $\Ceff = \langle E_1, E_2, H-E_1-E_2 \rangle$ and
$-K_2 = 3H-E_1-E_2$. Therefore, $\alpha(S_2)$ is the volume of
\begin{equation*}
  \begin{split}
    &\{(a_0,a_1,a_2) \in \RR^3 \mid 3a_0-a_1-a_2=1, a_1\ge 0, a_2 \ge
    0, a_0-a_1-a_2 \ge 0\}\\
    =&\{(a_0,a_1) \in \RR^2 \mid a_1 \ge 0, 3a_0-a_1-1 \ge 0, -2a_0+1 \ge 0\}\\
    =&\text{convex hull of $(1/3,0)$, $(1/2,0)$, $(1/2,1/2)$},
  \end{split}
\end{equation*}
which is a rectangular triangle whose legs have length $1/6$ and
$1/2$. Hence, $\alpha(S_2) = 1/24$, while $\alpha(S_1) = 1/6$,
$\alpha(\Pone \times \Pone)=1/4$, and $\alpha(\Ptwo)=1$ can also be
calculated directly, which completes the proof of
Theorem~\ref{thm:main}.
\end{proof}

\begin{remark}
  By the proof of \cite[Lemme~9.4.2]{MR1340296}, $\alpha(S_1) = 1/6$, and
  by the proof of \cite[Lemme~10.4.2]{MR1340296}, \[\alpha(S_2) = 1/3
  \cdot \vol\{(x_1,x_2) \in \RR^2_{>0} \mid x_1+x_2 \le 1/2\},\] which
  is clearly $\alpha_2 = 1/24$ and therefore agrees with our result.
  Note that the recursion formula does not hold for $r = 2$:
  \[\alpha(S_2) = \frac 1 {24} \ne \frac{N_2 \cdot \alpha(S_1)}{2\cdot (9-2)} =
  \frac 1{28}.\] The value $\alpha(S_4) = 1/(6\cdot 4!)$ was
  previously calculated in \cite[Section~1.3]{MR1909606}.
\end{remark}

\section{Singular Del Pezzo surfaces}\label{sec:singular}

For the classification of singular Del Pezzo surfaces, see
\cite{MR80f:14021}, \cite{MR89f:11083}, \cite{MR2227002}. It
turns out that in each degree, the surfaces can be divided into
different types according to the number and types of their
singularities and the number of lines.

For each type, we might have more than one isomorphism class (e.g.,
two for type $\Dfour$ in degree 3, and an infinite family for type
$\Aone$ in degree 3). However, the degrees in $\Pic(\tS_r)$ of the
negative curves, which generate the effective cone, and of $-K_r$ only
depend on the type. Therefore, $\alpha(S_r)$ depends only on the type.

More precisely, all this information is encoded for each type in its
extended Dynkin diagram of negative curves. For degree $\ge 4$, these diagrams
can be found in \cite[Propositions~6.1,~8.1,~8.3,~8.4]{MR89f:11083}. In
\cite{MR2227002}, the information for each degree is
encoded in a smaller number of diagrams. Also see
\cite{math.AG/0604194} and \cite[Chapter~6]{thesis} for detailed
information on several singular types of various degrees.

From the extended Dynkin diagram, we can derive a basis of
$\Pic(\tS_r)$, and $-K_r$ and all effective divisors in terms of this
basis as explained in \cite[Section~3]{math.AG/0604194}.

With this information, it is straightforward to calculate
$\alpha(S_r)$ for minimal desingularizations of all split types of
degree $\ge 3$. In practice, this task is significantly simplified by
the use of software such as $\polymake$.

For the Hirzebruch surface $F_2$, which is the unique singular Del Pezzo
surface of degree 8, we have $\alpha(F_2) = 1/8$.

For $S_r$ of types $\Dfour$ and $\Dfive$ of degree 4 and $\Esix$ of
degree 3, the constant $\alpha(S_r)$ has been calculated while proving
Manin's conjecture for these surfaces; see \cite{math.NT/0604193},
\cite{math.NT/0412086}, \cite{math.NT/0509370}, respectively. Note
that in these three cases, the calculation of $\alpha(S_r)$ is
particularly simple as the effective cone is \emph{simplicial}, i.e.,
the number of generators of $\Ceff$ equals the rank of $\Pic(\tS_r)$.

\begin{table}[ht]
  \centering
\[\begin{array}{|c||c|c|c|c|}
  \hline
  \text{type} & \text{singularities} 
  & \text{lines} & \text{generators of $\Cdual$} & \alpha\\
  \hline
  \hline
  0 & - & 3 & 3 & 1/24\\
  i & \Aone & 2 & 3 & 1/48\\
  \hline
\end{array}\]
  \caption{Del Pezzo surfaces of degree 7}
  \label{tab:degree_7}
\end{table}

\begin{table}[ht]
  \centering
\[\begin{array}{|c||c|c|c|c|}
  \hline
  \text{type} & \text{singularities} 
  & \text{lines} & \text{generators of $\Cdual$} & \alpha\\
  \hline
  \hline
  0 & - & 6 & 5 & 1/72\\
  i & \Aone & 4 & 5 & 1/144\\
  ii & \Aone & 3 & 4 & 1/144\\
  iii & 2\Aone & 2 & 4 & 1/288\\
  iv & \Atwo & 2 & 4 & 1/432\\
  v & \Atwo+\Aone & 1 & 4 & 1/864\\
  \hline
\end{array}\]
  \caption{Del Pezzo surfaces of degree 6}
  \label{tab:degree_6}
\end{table}

\begin{table}[ht]
  \centering
\[\begin{array}{|c||c|c|c|c|}
  \hline
  \text{type} & \text{singularities} 
  & \text{lines} & \text{generators of $\Cdual$} & \alpha\\
  \hline
  \hline
  0 & - & 10 & 10 & 1/144\\
  i & \Aone & 7 & 9 & 1/288\\
  ii & 2\Aone & 5 & 8 & 1/576\\
  iii & \Atwo & 4 & 7 & 1/864\\
  iv & \Atwo+\Aone & 3 & 7 & 1/1728\\
  v & \Athree & 2 & 5 & 1/3456\\
  vi & \Afour & 1 & 5 & 1/17280\\
  \hline
\end{array}\]
  \caption{Del Pezzo surfaces of degree 5}
  \label{tab:degree_5}
\end{table}

\begin{table}[ht]
  \centering
\[\begin{array}{|c||c|c|c|c|c|}
  \hline
  \text{type} & \text{singularities} 
  & \text{lines} & \text{generators of $\Cdual$} & \alpha\\
  \hline
  \hline
  0 & - & 16 & 26 &  1/180\\
  i & \Aone & 12 & 22 & 1/360\\
  ii & 2\Aone & 9 & 19 & 1/720\\
  iii & 2\Aone & 8 & 17 & 1/720\\
  iv & \Atwo & 8 & 16 & 1/1080\\
  v & 3\Aone & 6 & 15 & 1/1440\\
  vi & \Atwo+\Aone & 6 & 15 & 1/2160\\
  vii & \Athree & 5 & 11 & 1/4320\\
  viii & \Athree & 4 & 10 & 1/4320\\
  ix & 4\Aone & 4 & 4 & 1/2880\\
  x & \Atwo+2\Aone & 4 & 13 & 1/4320\\
  xi & \Athree+\Aone & 3 & 10 & 1/8640\\
  xii & \Afour & 3 & 9 & 1/21600\\
  xiii & \Dfour & 2 & 6 & 1/34560\\
  xiv & \Athree + 2\Aone & 2 & 10 & 1/17280\\
  xv & \Dfive & 1 & 6 & 1/345600\\
  \hline
\end{array}\]
  \caption{Del Pezzo surfaces of degree 4}
  \label{tab:degree_4}
\end{table}

\begin{table}[ht]
  \centering
\[\begin{array}{|c||c|c|c|c|c|}
  \hline
  \text{type} & \text{singularities} 
  & \text{lines} & \text{generators of $\Cdual$} & \alpha\\
  \hline
  \hline
  0 & - & 27 & 99 & 1/120\\
  i & \Aone & 21 & 78 & 1/240\\
  ii & 2\Aone & 16 & 62 & 1/480\\
  iii & \Atwo & 15 & 52 & 1/720\\
  iv & 3\Aone & 12 & 49 & 1/960\\
  v & \Atwo+\Aone & 11 & 43 & 1/1440\\
  vi & \Athree & 10 & 32 & 1/2880\\
  vii & 4\Aone & 9 & 38 & 1/1920\\
  viii & \Atwo + 2\Aone & 8 & 35 & 1/2880\\
  ix & \Athree + \Aone & 7 & 27 & 1/5760\\
  x & 2\Atwo & 7 & 31 & 1/4320\\
  xi & \Afour & 6 & 21 & 1/14400\\
  xii & \Dfour & 6 & 16 & 1/23040\\
  xiii & \Athree + 2\Aone & 5 & 23 & 1/11520\\
  xiv & 2\Atwo+\Aone & 5 & 18 & 1/28800\\
  xv & \Afour+\Aone & 4 & 18 & 1/28800\\
  xvi & \Afive & 3 & 13 & 1/86400\\
  xvii & \Dfive & 3 & 11 & 1/230400\\
  xviii & 3\Atwo & 3 & 21 & 1/25920\\
  xix & \Afive+\Aone & 2 & 13 & 1/172800\\
  xx & \Esix & 1 & 7 & 1/6220800\\
  \hline
\end{array}\]
  \caption{Del Pezzo surfaces of degree 3}
  \label{tab:degree_3}
\end{table}

\begin{remark}
  We have two different types with the same singularities in three
  cases ($\Aone$ in degree 6, $2\Aone$ and $\Athree$ in degree 4).
  The two types in each pair can be distinguished by their number of
  lines. Therefore, the two different types have different effective
  cones. However, both types in each pair have the same constant
  $\alpha$. It is unclear whether this is more than a coincidence.
\end{remark}

\begin{remark}
  For singular Del Pezzo surfaces $S_r$ of degree 2 and 1, we do not
  calculate $\alpha(S_r)$ since the number of different types is much
  larger than in degree $\ge 3$.
  
  For surfaces of degree 1 whose effective cone has many generators,
  the triangulation of the nef cone might be too complicated for a
  computation of $\alpha$ using $\polymake$. We expect this to happen
  in case of ``mild'' singularities (e.g., of type $\Aone$, $\Atwo$ or
  $2\Aone$). Here, the nef cone should be almost as complicated as in
  the smooth case, and some help ``by hand'' might be necessary.
\end{remark}

\bibliographystyle{alpha}

\bibliography{alpha}

\end{document}